\documentclass{article}

\usepackage{amsfonts}

\def\*{{\ast}}

\def\comp{\circ}

\newtheorem{proposition}{Proposition}
\newtheorem{theorem}[proposition]{Theorem}
\newtheorem{lemma}[proposition]{Lemma}
\newtheorem{corollary}[proposition]{Corollary}
\newtheorem{definition}[proposition]{Definition}

\newcommand{\evalAt}{\,\vrule width.2pt\relax\,{\vphantom \sum}}

\newcommand{\eqv}{\simeq}
\newcommand{\exchg}{\equiv}

\newcommand{\umbDot}{\!\mathbin{\raisebox{-1pt}{$\scriptscriptstyle\bullet$}}\!}
\renewcommand{\.}{\mathchoice{\umbDot}{\umbDot}{\,\umbDot\,}{\,\umbDot\,}}

\newcommand{\E}{\mathop{\hbox{\bf E}}}
\newcommand{\multiplicative}{multiplicative }

\newcommand{\prf}{\par{\noindent\em Proof.}}
\newenvironment{proof}{\prf \setlength{\parindent}{3ex}}{\qed\par\noindent}
\def\squarebox#1{\hbox to #1{\hfill\vbox to #1{\vfill}}}
\newcommand{\qedbox}{\vbox{\hrule\hbox{\vrule\squarebox{.667em}\vrule}\hrule}}
\newcommand{\qed}{\nopagebreak\mbox{}\hfill\qedbox\smallskip}

\newcommand{\Q}{{\mathbb Q}}
\newcommand{\F}{{\bf F}}
\renewcommand{\k}{{\bf k}}
\newcommand{\A}{{\cal A}}
\newcommand{\C}{{\mathbb C}}

\author{Brian D. Taylor\\
Wayne State University\\
Detroit, MI 48202, USA}
\title{Umbral presentations for polynomial sequences.}
\date{April 30, 1999}

\begin{document}



\maketitle

\begin{abstract}
Using random variables as motivation,
this paper presents an exposition of formalisms developed in~\cite{RT1,RT2}
for the classical umbral calculus. A variety of examples are presented,
culminating in several descriptions of sequences of binomial type
in terms of umbral polynomials.
\end{abstract}

\section{Introduction}

The system of calculation now known as the ``umbral calculus''
originated with Blissard in the nineteenth century in informal calculations
involving the ``lowering'' and ``raising'' of exponents.
The work of Rota and his collaborators in~\cite{MR,RKO,RomRot} and other works
formalized these methods in the modern language of linear operators and
Hopf algebras.
While this clarified the underlying theory, it rendered the original
nineteenth century work no more easy to read or check. In~\cite{RT1,RT2,T},
the original classical notation was revived and extended so as to be rigorous by
modern standards. In~\cite{RST}, the first attempts were made to apply
this newly revived classical umbral calculus to one of the most significant
successes of the modern theory, namely the study of sequences of binomial type.

The purposes of this paper are twofold. Since much of this paper is expository,
no prior knowledge of umbral calculus in any of its guises is assumed.
To start, we develop the modern formulation of the classical umbral calculus in 
analogy with the idea of a random variable. This renders the definitions of~\cite{RT1,RT2}
transparent.   
In the second part of this paper, we introduce a new operation on umbrae
arising naturally from the analogy to random variables.
We show that all sequences of binomial type and all umbral maps arise directly
from the application of this operation.
We further apply the tools of classical umbral calculus developed in the
first part of the paper to provide several other compact presentations for sequences of 
binomial type.

The author has attempted to document at least the recent history of the main results 
and definitions contained herein. 

\section{Random variables and the classical umbral calculus}

Fundamental to the classical umbral calculus is the idea of associating
a sequences of numbers $a_0,a_1,a_2,\ldots$ to an ``umbral variable'' $\alpha$ which
is said to represent the sequence. To be slightly more formal, the umbral
calculus relies on associating the sequence $a_0,a_1,a_2,\ldots$ to the sequence
$1,\alpha,\alpha^2,\alpha^3,\ldots$ of powers of $\alpha$.
 
This kind of association is 
familiar in modern mathematics: To any random variable $G$, we associate
a sequence of numbers $1,g_1,g_2,\ldots$ where $g_i$ is the $i$th moment
of $G$. Specifically, we are defining the sequence
$1,g_1,g_2,\ldots$ by applying the expectation operator, $\E$,
componentwise to the sequence $1,G,G^2,\ldots$ consisting of powers of 
the random variable $G$.

We will proceed to carry this analogy further, calculating with random
variables in precisely the way we will later be using umbral variables.  
We start by letting $G$ be a random variable distributed uniformly over
the interval $[0,1]$. The sequence of moments $1,g_1,g_2,\ldots$
associated to $G$ is thus given by $g_n=\int_0^1 g^n\,dg={1\over n+1}$.

If we let $p(t)\in\C[t]$ be any polynomial with complex coefficients,
it is immediate that 
\begin{equation}\label{firsteval}
 \E[p'(G)] =\int_0^1 p'(g)\,dg = \Delta p(0) \, ,
\end{equation}
where $\Delta$ is the forward difference operator $\Delta p(t) =p(t+1)-p(t)$.
A comment on notation: since $\Delta$ is defined as an operator on the ring of 
polynomials in $t$, $\Delta p(0)$ can only be interpreted as $\Big(\Delta p\Big)(0)$
or as $0$. We adopt the former reading.
Since the calculation in equation~\ref{firsteval} only required that $p(t)$ was 
differentiable, it could just as well have been carried out for a polynomial with
coefficients in some other integral domain or indeed for a polynomial in $t$
whose coefficients contained various random variables {\em which were independent}
of $G$. So suppose $G'$ is a random variable independent of and identically
distributed to $G$. Consider $p(G+G')$ as a polynomial $q(G)$ in $G$ with 
coefficients in $\C[G']$. Using equation~\ref{firsteval}, the expected value, 
averaging over values of $G$, of $p''(G+G') = q''(G)$ is $\Delta q(0)=\Delta p'(G')$.
Applying  equation~\ref{firsteval} again, recalling that the derivative $D$ and $\Delta=e^D-I$
commute and that $G'$ is identically
distributed to $G$, gives $\E[\Delta p'(G')]= \Delta^2 p(0)$. Somewhat more
suggestively, this calculation can be written
\begin{equation}\label{secndeval}
\E[p''(G+G')] = \E[\Delta p'(G')] = \E[\Delta p(0)] \, .
\end{equation}

The first property to observe here is that the independence of $G$ and $G'$ matters,
as  
\[ \E[p''(G+G)] = \E[\Delta p''(2G)] = {p'(2\cdot1)-p'(0)\over 2} \, ,\]
since $p''(2t)=D_t\, {p'(2t)\over 2}$. where $D_t$ is the derivative with respect to $t$. 

The second important property is the rather trivial observation that we can calculate
the expectation of a polynomial in several random variables all independent and 
identically distributed to $G$ simply by knowing the moments of $G$. 
For example,
by equation~\ref{secndeval}, we know that $\E[20(G+G')^3]=\Delta^2 t^5\evalAt_{t=0}$.
This could be evaluated directly. Alternatively, and denoting 
$E[G^k]$ by $G_k$, it could be evaluated as
\begin{eqnarray*}
20\cdot\E[G^3+3G^2G'+3GG'^2+G'^3 ] &=& 20(G_3+3G_2G_1+3G_1G_2+G_3 )\\
&=& 20(2G_3+6G_2G_1)\\
&=&20\left({2\over4}+{6\over 3\cdot2}\right)\,.\\
\end{eqnarray*}

Applying the above observations quickly reconstructs the moment generating 
function $E[e^{Gz}]$ for $G$. Since $D_t\, {e^{tz}\over z}= e^{tz}$,
equation~\ref{firsteval} implies $E[e^{Gz}]={e^z-e^0\over z}$.

For the duration of the next calculation, we are going to make some assumptions
that simply do not hold within the confines of probability theory. The remainder
of this section will be devoted to describing how to replace random variables
with ``umbral variables'' in a way that makes the following calculations legitimate.
So, for the moment, let us assume that there is an object $B$ that behaves much like
a random variable. Let us call this object an ''umbral variable.'' We treat it just
like a random variable, but stipulate both 
that it is independent of $B$ and that $B+G=0$.
With these stipulations, we find that by independence
\[ \E[e^{(G+B)z}=\E[e^{Gz}e^{Bz}]=\E[e^{Gz}]\E[e^{Bz}]={e^z-1\over z} \E[e^{Bz}] \, .\]
But since $G+B=0$, the left-hand side above is just $1$. We have just calculated
that $\E[e^{Bz}]={t\over e^t-1}$. Since this is the exponential generating function,
$\sum_{k\ge0} B_k {z^k\over k!}$,
for the Bernoulli numbers, $B_k$, we find that if the above calculation can be
made rigorous, then $\E[B^k]=B_k$.

Since calculations such as the previous are too useful to abandon
(see~\cite{RT2} for a variety of examples involving the Bernoulli numbers)
we define {\em umbral variables} or {\em umbrae} which formalize
the roles of both the random variable $G$ and the new object $B$ in the
preceding calculation.

Just as a random variable is usually capitalized, we will typically
distinguish our umbrae by writing them as Greek letters, e.g. $\alpha,\beta,
\alpha''',\ldots$. Let us denote the collection of whichever umbral variables
we will be using by $\A$. See~\cite{RT1,RT2} for the relevant, and straightforward,
technical details. 
In practice, when we introduce a new umbra, say $\alpha$, we specify explicitly or 
implicitly how $\E$ acts on it, namely what values $\E[\alpha^k]$ takes for each $k$.
Any two distinct umbrae in $\A$, say $\alpha$ and $\gamma$ or $\alpha$ and $\alpha'$
will acts as independent random variables, regardless of how $\E$ acts on them.
Generalizing, any collection of distinct umbrae will behave as independent random
variables.

Formally, this can be accomplished by defining a linear {\em evaluation} map, 
$\E:\F[\A]\rightarrow\F$ on the polynomials $\F[\A]$ in the umbrae with coefficients
in a suitably chosen commutative ring $\F$. We require that $\E$ is $\F$-linear,
$\E[1]=1$, and  $\E[M\cdot M']=\E[M]\E[M']$ for any two monomials, $M$ and $M'$ in
$\C A$ such that no umbra appears to nonzero power in both $M$ and $M'$. This 
map was called ${\bf eval}$ in~\cite{RT1,RT2}. 

We call $p,q\in\F[\A]$ {\em umbrally equivalent}, written $p\eqv q$
when $\E[p]=\E[q]$. Analogous to the notion of identically distributed
random variables, we define $p,q\in\F[\A]$ to be {\em umbrally exchangeable}, 
when $p^k\eqv q^k$ for all $k\ge0$. If for example we have $\alpha\exchg 3$,
then $\E[\alpha^k]=3^k$; this is consistent, in the analogy to random variables,
with considering  $\alpha$ analogous to a random variable which always takes
on the value $3$. We not that equality implies exchangeability which implies
umbral equivalence. The converses are false.

We define  $p,q\in\F[\A]$ to be 
{\em independent} when no umbra appears in both $p$ and $q$. More formally,
an umbra that appears to a nonzero power in some monomial with nonzero coefficient
in $p$ does not appear to a nonzero power in any monomial with nonzero coefficient
in $q$.
For example,  $\alpha^2+\alpha \alpha'$ and $\beta\beta'^2-\beta+\alpha'''$
are independent, but $\alpha^2+\alpha \alpha'$ and $\beta\beta'^2-\beta+\alpha$
are not independent. Nor are the falling factorials $(\alpha)_{(n)}$
and $(\alpha)_{(n-1)}$ independent. On the other hand,
$\alpha t^2+\beta$ and $\gamma-t$ are independent where $t$ is some element of $\F$.

In random variables, we can usually substitute one identically
distributed random variable for another (modulo independence constraints).
Similarly, we can substitute exchangeable umbrae as per the following lemma.
\begin{lemma}[RT1]\label{subsLemma}
If a polynomial $p(t)$ is independent of two
exchangeable umbrae
$\alpha$ and $\alpha'$, then $p(\alpha)\exchg p(\alpha')$.
\qed\end{lemma}
This substitution lemma holds equally well if $\alpha$ or $\alpha'$ is replaced by an 
umbral polynomial $p\in\F[\A]$.

For the duration of this paper, we let $\varepsilon$ be an umbra such that
$\varepsilon^k\eqv\delta_{0,k}$ where $\delta$ in the Kronecker delta. As
long as we work with polynomials in the umbrae, there is no harm in defining
$0^0=1$. Under this convention we consider $0\exchg\varepsilon$. This is consistent
with the convention that $1$ is the $0$th moment of a random variable which always takes
on value $0$.

To pick up our earlier example, let $\gamma$ be an umbra such that
$\gamma^i\eqv {1\over n+1}$ and let $\beta$ be an umbra such that
$\gamma+\beta\exchg0$. It is an easy exercise to see that given any
umbra $\gamma$ such an {\em inverse umbra} can be found recursively.
Here we have made formal in $\gamma$ and $\beta$
exactly the properties we had assumed for $G$ and $B$.
We have $\gamma$ and $\beta$ are independent and $\beta+\gamma\exchg0$.

Extending our notions of independence, exchangeability, and equivalence
coefficientwise to formal power series (see~\cite{T} for a general but technical
treatment) we can duplicate the computation we did for $G$ and $B$.
By the substitution lemma we have $e^{(\beta+\gamma)z}\eqv e^{0\cdot z}=1$.
Thus we have $e^{\beta z}e^{\gamma z}\eqv 1$. By independence, this implies
that $e^{\beta z}{e^z-1\over z}\eqv 1$ and hence (by linearity) that 
$e^{\beta z}\eqv{z\over e^z-1}$. From a technical viewpoint, there are
a number of ways to justify the
first step in the preceding sentence. 
The most direct solution is to apply the substitution lemma coefficientwise to the formal
power series in $z$.
A general approach which views multiplication
by $e^{\gamma z}{z\over e^z-1}$ as a linear operator ``equivalent'' to the identity
is given in~\cite{T}. The intuition behind both these proofs is that since 
$e^{\gamma z}$ and $e^{\beta z}$ are independent, $\E$ can be 
applied in two stages first to $\gamma$ and then to $\beta$,
analogously to finding the expectation by 
first averaging over one random variable and then over another
independent random variable.

As a demonstration of these techniques we rewrite in modern umbral notation the first example
in~\cite{Bl}, one in the series of papers
in which Blissard during the 1860's introduced his ``representative notation''---the umbral
calculus. To point out just how closely the modern language captures Blissard's
19th century original, we present most of this example in Blissard's own words.

\begin{minipage}[t]{4.25in}
Blissard starts with the problem
``Required to expand $\left\{x\over\log(1+x)\right\}^m$.''
He then lets
 ``$\left\{x\over\log(1+x)\right\}^m=1+P_1x+P_2x^2+\cdots+P_nx^n+\&c.$''
and defines $U_n$ to be the coefficient of $\theta^n\over n!$ in 
$\left(e^\theta-1\over\theta\right)^m$ where $\theta$ is an ordinary variable. 
He observes that
\begin{eqnarray*}
``\quad\left(e^\theta-1\over\theta\right)^m &=&
1+U_1\theta+U_2{\theta^2\over1\cdot2}+\cdots+U_n{\theta^n\over1\cdot2\cdots n}+\&c.\\
&=&e^{U\theta} \hbox{(by representative notation)\,.\quad''}
\end{eqnarray*}
In modern language, he is letting $U$ be an umbra such that $U^n\eqv U_n$,
also his ``$=$'' would be replaced with ``$\eqv$''.
The next operation takes place purely on the level of formal power series.
Blissard substitutes $\log(1+x)$ for $\theta$ and finds that
\[ ``\quad 
\left\{x\over\log(1+x)\right\}^m=(1+x)^U
\quad\hbox{''}\]
where again the only change necessary to modernize his work is to replace
 ``$=$'' with ``$\eqv$''. If we
``equate coefficients of $x^n$, then 
$P_n={U(U-1)(U-2)\cdots(U-n+1)\over 1\cdot2\cdot3\cdots n}$'';
again we would replace ``$=$'' with ``$\eqv$''.
\end{minipage}
\vskip.1in

The preceding formula for $P_n$ has the advantage of being extremely compact.
Blissard concludes with an expansion of it, and we shall proceed likewise, though
our precise techniques are somewhat more umbral than those Blissard used.

With $\gamma$ as before, we have ${e^\theta-1\over \theta}\eqv e^{\gamma\theta}$.
Thus $\left({e^\theta-1\over \theta}\right)^m\eqv e^{(\gamma'+\gamma''+\cdots+\gamma''')\theta}$
where $\gamma'+\gamma''+\cdots+\gamma'''$ is a sum of $m$ distinct (and thus independent)
umbrae each exchangeable with $\gamma$. We conclude that 
$U\exchg \gamma'+\gamma''+\cdots+\gamma'''$ and thus 
$P_n\eqv {\gamma'+\gamma''+\cdots+\gamma'''\choose n}$.
We conclude with the following formula for evaluating the powers $U^n$.
Since $D_t^m \,t^{m+n}=(m+n)_{(m)} t^n$, we have, generalizing calculation~\ref{secndeval},
that 
\[ (m+n)_{(m)} U^n\eqv (m+n)_{(m)} (\gamma'+\gamma''+\cdots+\gamma''')^n \eqv
\Delta^m\,0^{m+n}\,. \]
This last is better known as $m!\,S(m+n,m)$, where $S(n,k)$ is the 
Stirling number of the second kind counting the number of 
set partitions of an $n$-set into $k$ parts.
So $U^n\eqv S(m+n,m)/{m+n\choose m}$.
We can use this, together with the expansion of the falling factorials
in terms of Stirling numbers of the first kind to derive
\[{\gamma'+\cdots+\gamma'''\choose n}=
{1\over n!}\sum_{k=0}^n s(n,k)(\gamma'+\cdots+\gamma''')^k
\eqv
{1\over n!}\sum_{k=0}^n {s(n,k)  S(m+k,m)\over {m+k\choose m}} \, .
\]

\section{Umbral presentations of Appell sequences}

Historically, the objects of interest in umbral computations were of course
sequences of numbers or polynomials. 
For our present purposes, this means that we will primarily be
studying the ``moments'' $E[\alpha^k]$ of an umbra $\alpha$.
We say that the umbra $\alpha$ {\em represents} a sequence $a_0,a_1,a_2,\ldots$,
$a_i\in\F$, when $\alpha^k\eqv a_k$ for all integers $k\ge0$. Necessarily this
implies that $a_0=1$. An {\em umbral presentation} of a sequence 
$a_1,a_2,\ldots$ of elements in $\F$ is any sequence $q_1,q_2,\ldots$ of polynomials in
$\F[\A]$ such that $q_i\eqv a_i$ for $i\ge0$.
Throughout this paper we freely assume that, for any sequence in $\F$, we can find 
infinitely many umbrae representing the given sequence.

Now let $\F$ be $\k[x,y]$ where $\k$ is a commutative ring containing $\Q$.
The remainder  of this paper will focus on umbral presentations for sequences 
of polynomials. For example for any umbra $\alpha$, we can define a sequence of 
polynomials $s_n(x)$, $n=0,1,2,\ldots$, by $s_n(x)\eqv (x+\alpha)^n$. This
definition immediately yields the calculation
\begin{equation}\label{appellcalc}
s_{n}(y+x) \eqv (y+x+\alpha)^n = \sum_i {n\choose i} y^i (x+\alpha)^{n-i} \eqv 
\sum_i {n\choose i} y^i s_{n-i}(x) \, .
\end{equation}

A sequence of polynomials $s_0(x),s_1(x),s_2(x),\ldots$ is said to be an
{\em Appell sequence} when it satisfies the identity, 
\begin{equation}\label{appelldef} 
s_{n}(y+x) = \sum_i {n\choose i} y^i s_{n-i}(x)
\end{equation}
given by equation~\ref{appellcalc}
for all $n\ge0$. We shall call an Appell sequence $s_n(x)$ {\em normalized}
when $s_1(x)$ is monic. Any Appell sequence may be rewritten as a normalized
Appell sequence by replacing $s_n(x)$ with $s_n(x)/s_1'(0)$. Here, $s_1'(x)$
is the first derivative of $s(x)$.  We hold with this notation for derivatives
throughout this paper.  In the literature, Appell sequences
are frequently defined to be normalized.

\begin{proposition}[\cite{RT2}]
A sequence, $s_0(x),s_1(x),\ldots$, of polynomials in $\k[x]$ with $s_n(x)$ having
degree $n$ is a normalized Appell sequence iff there exists an umbra $\alpha$ such that
$ s_n\eqv (x+\alpha)^n$ 
for $n\ge0$.
\end{proposition}
\begin{proof}
The if direction is given by calculation~\ref{appellcalc}. 

(only if): Replacing $x$ with $0$ in the defining equation~\ref{appelldef} shows
that in an Appell sequence, each polynomial $s_n(y)$ can be recovered
from the sequence of values $s_0(0),s_1(0),\ldots$.
Choosing an umbra $\alpha$ that represents this sequence guarantees
$(x+\alpha)^n\eqv s_n(x)$.
\end{proof}

Similarly we have the standard result that a sequence $s_0(x),s_1(x),\ldots$ of polynomials,
$s_n(x)$ of degree $n$, is an Appell sequence iff 
\begin{equation}\label{diffAppell}
s_n'(x)=n\cdot s_{n-1}(x) 
\end{equation}
for all $n\ge0$. The only if direction follows since
$D_x (x+\alpha)^n = n(x+\alpha)^{n-1}$. To show the if direction, we observe 
that any sequence of polynomials satisfying equation~\ref{diffAppell} is determined by
the sequence of values of $s_0(0), s_1(0), s_2(0),\ldots$ and apply
the argument in the preceding proof. 

The sequences of polynomials with which we will most concerned in this
paper are those of ``binomial type,'' i.e. sequences of polynomials which
satisfy an analog of the binomial theorem. Before approaching this topic
however, we need to lift another tool from random variables to umbral
calculus.

\section{Sums of random variables and the ``dot'' operation on umbrae}
 
Suppose that $X$ is a random variable. If $n$ is a positive integer,
one can of course run $n$ trials of $X$ and sum the results. Denote the
sum by a new random variable $n\.X$. Thus $n\.X$ has the same distribution
as $X_1+X_2+\cdots+X_n$ where the $X_i$ are all independent and identically
distributed to $X$. In~\cite{RT1,RT2} the corresponding notion $n.\alpha$
was defined for an arbitrary umbra $\alpha$. In particular $n\.\alpha$
is itself an umbra and it is defined to be exchangeable with
$\alpha_1+\alpha_2+\cdots+\alpha_n$ where $\alpha_i\exchg \alpha$ for
each $i$. Similarly, for any umbral polynomial $p\in\F[\A]$ we define
a new umbra $n\.p$  which is exchangeable with $n\.\gamma$ where $\gamma$
is any umbra satisfying $\gamma\exchg p$.
It is worth emphasizing that $n\.p$ is itself an umbra.  
Thus $\alpha$ $5\.\alpha$, $3\.\alpha$, and $5\.(\alpha+\beta)$ are
all distinct (hence all independent.) It is however clear from the
definitions that $5\.(\alpha+\beta)\exchg 5\.\alpha + 5\.\beta$.

We recall a technical consideration from~\cite{RT1,RT2}.
The set of all umbrae $\A$ will be decomposed as a disjoint
union $\A=\A_0\uplus \A_1$, umbrae in $\A_1$ are called
{\em auxiliary} umbrae. Umbrae of the form $n\.p$ are
auxiliary umbrae. This detail will be given more attention below.

It is an easy observation that if $g(z)\eqv e^{\alpha\, z}$ (this object
is analogous to the moment generating function of a random variable)
then 
\[ e^{n\.\alpha\, z}\exchg e^{\alpha_1+\cdots+\alpha_n\, z} =
   \prod_{i=1}^n e^{\alpha_i\, z} \eqv g(z)^n \, ,\]
where the last equivalence uses the independence of the $e^{\alpha_i z}$'s.

Even more directly, we see that 
\begin{equation}\label{bintake1}
(m+n)\.\alpha \exchg m\.\alpha + n\.\alpha'
\end{equation}
where $\alpha'\exchg\alpha$. As a consequence, if for each positive integer
$n$ we define a sequence $f_0(n),f_1(n),f_2(n),\ldots$, by
$f_i(n)\eqv (n\.\alpha)^i$, then equation~\ref{bintake1}
implies that 
\begin{displaymath}
 ((m+n)\.\alpha)^k \eqv (m\.\alpha + n\.\alpha')^k = 
        \sum_i {k\choose i} (m\.\alpha)^i (n\.\alpha')^{k-i}
\end{displaymath}
hence
\begin{equation}\label{bintake1exp}
f_k(m+n) = \sum_i {k\choose i} f_i(m)\, f_{k-i}(n) \,\, .
\end{equation}
This kind of generalized binomial theorem will be explored further
in the next section.

By way of introduction to the first new definition of this paper, we consider
the following generalization of $n\.X$. Let $X$ be some random variable and let $Y$ be
a random variable which only takes positive integer values. Run one trial
of $Y$, then run $Y$ trials of $X$ and sum the results. 
We define $Y\.X$ to be a new random variable whose distribution is identical
to $X_1+X_2+\cdots+X_Y$; 
for convenience, we are defining the  $X_i$'s to be independent random variables 
identically distributed to $X$.
Observe that if $Z$ is another random variable
taking only takes positive integer values, then according the preceding definition
$(Y+Z)\.X$ has the same distribution as $X_1+X_2+\cdots+X_Y+X_{Y+1}+\cdots+X_{Y+Z}$
and hence as $Y\.X_1 + Z\.X_2$. On the other hand, $X\.(Y+Z)$ does not in general
have the same distribution as $X_1\.Y+X_2\.Z$; in the first expression only one
trial of $X$ is made and in the second one makes two trials of $X$. Of course
if $X$ always returns the same value, say $X=n$, this causes no trouble and
$n\.(X+Y)$ is identically distributed to $n\.X + n\.Y$.

Just as the definition $n\.X$  for random variables
generalizes to $n\.\gamma$ for umbrae, we would like 
a generalization of the random variable $X\.Y$ to umbrae.
This generalization should satisfy results analogous to those observed above for
random variables.

The generalization relies on a simple observation first applied to the umbral
calculus by Nigel Ray in~\cite{Ra2}.
\begin{proposition}\label{Ray-poly}
If $\gamma\in\A_0$ is an umbra, and $n$ is a positive integer,
then \goodbreak
$E[(n\.\gamma)^k]$ is a polynomial in $n$.
\end{proposition}
\begin{proof}
This is equivalent to the observation that if 
$g(z)$ is in $\k[[z]]$, the ring of formal power series in $z$, and if
$g(z)\eqv e^{\gamma z}$, 
then  $e^{n\.\gamma \, z} \eqv g(z)^n = e^{n\log( g(z) )}$ and
the coefficient of $z^k/k!$ in the last is a polynomial in $n$.
\end{proof}

Alternately, we could have observed that if $\gamma^i\eqv a_i$,
and each $\gamma_i\exchg \gamma$ for $i=1\ldots n$
then $a_1^{i_1}a_2^{i_2}\cdots a_k^{i_k}$ appears in $\E[(\gamma_1+\cdots+\gamma_n)^k]$
as many times as there are
monomials in the expansion of $(\gamma_1+\cdots+\gamma_n)^k$
containing exactly $i_j$ $j$th powers.
But this says that $(n\.\gamma)^k$ is umbrally equivalent to
\begin{equation}\label{ProbCoeffs}
(\gamma_1+\cdots+\gamma_n)^k \eqv
\sum_{i_1,\ldots,i_k}
{n\choose i_0,\ldots,i_k}
{{\displaystyle k}\choose 
    \underbrace{1,\ldots,1}_{\hbox{$i_1$ times}},\ \ldots,\  
    \underbrace{k,\ldots,k}_{\hbox{$i_k$ times}}}
a_1^{i_1}\cdots a_k^{i_k}
\end{equation}
where $i_0=n-(i_1+\ldots+i_k)$ and 
${n\choose n-(i_1+\ldots+i_k),i_1,\ldots,i_k}$ 
is a polynomial of degree $(i_1+\ldots+i_k)$ in $n$.

But this is precisely what we need to make sense of replacing $n$ with $\alpha$.
\begin{definition}\label{umbraldot}
Let $\alpha,\gamma\in\A_0$ be umbrae, and define $g(z)\in\F[[z]]$ by
$g(z)\eqv e^{\gamma z}$. Let $q_{\gamma,k}(n)$ be the coefficient
of  $z^k/k!$ in $e^{n\log( g(z) )}$.  Define $\alpha\.\gamma\in\A_1$ to
be a new auxiliary umbra,
such that $(\alpha\.\gamma)^k\eqv q_{\gamma,k}(\alpha)$.

\noindent
In general, if $p,q\in\F[\A_0]$ are umbral polynomials, we define
an auxiliary umbra $p\.q\in\A_1$ by $p\.q\exchg \alpha\.\beta$
where $\alpha\exchg p$ and $\beta\exchg q$.
\end{definition}
Equivalently we could have defined $(\alpha\.\gamma)$ by replacing $n$ with
$\alpha$ in equation~\ref{ProbCoeffs}.

The definition immediately implies that
$e^{\alpha\.\gamma\, z}\eqv e^{\alpha\log( g(z) )} = ( g(z) )^\alpha$.
Similarly, if $a(z)\in\k[[z]]$ is defined by $a(z)\eqv e^{\alpha z}$,
then $e^{\alpha\.\gamma\, z}\eqv a\big(\log(g(z))\big)$.
It is a straightforward exercise in probability theory to show
that if $a(z)$ is the moment generating function of a random variable
$A$ taking only positive integer values and if $g(z)$ is the moment
generating function of a random variable $B$ then 
$B_1+\cdots+B_A$ also has moment generating function $a\big(\log(g(z))\big)$.

It follows that under this definition $0\.\gamma\exchg \varepsilon\exchg 0$
which is what one would expect from the analogy to random variables.

We now state the promised analogues to the standard results
on random variables.
\begin{proposition}\label{leftLin}
Let $p,q,r\in\F[\A]$ be umbral polynomials.
If $p,q$ are independent then 
$(p+q)\.r \exchg p\.r + q\.r$.
\end{proposition}
\begin{proof}
By definition, and the substitution lemma (Lemma~\ref{subsLemma}),
it suffices to prove that 
$(\alpha+\beta)\.\gamma \exchg \alpha\.\gamma + \beta\.\gamma$
for any distinct umbrae $\alpha,\beta,\gamma$, 
i.e. that the $k$th powers of each side of the
displayed equation are umbrally equivalent for all $k\ge0$.
Letting $q_{\gamma,k}(n)$ be the polynomials from Definition~\ref{umbraldot},
it suffices to show, for all $k\ge0$, that the identity
$q_{\gamma,k}(\alpha+\beta)=
\sum_i {k\choose i} q_{\gamma,i}(\alpha) q_{\gamma,k-i}(\beta)$
holds purely on the level of polynomials in variables $\alpha,\gamma$.
But this follows since
equation~\ref{bintake1exp} says this identity holds with $\alpha,\gamma$
replaced by any pair of positive integers.
\end{proof}

As remarked above, we cannot expect that $p\.(q+r)\exchg p\.q + p\.r$
will hold in general. However, we record the special case where $p$
involves no umbrae.
\begin{proposition}\label{specialRightLin}
Let $a$ be an element of $\F$. Let $q,r\in\F[\A_1]$ be umbral polynomials.
If $q,r$ are independent, then $a\.(q+r)\exchg a\.q + a\.r$.
\end{proposition}
\begin{proof}
The result holds when $a$ is any integer. Repeating
the argument in the proof of Proposition~\ref{leftLin}
shows the identity holds when interpreted in terms of
polynomials in $a$.
\end{proof}
The importance of independence is illuminated if we examine
what fails on replacing $a$ with $\alpha$.
and trying to prove that $a\.(\beta+\gamma)\exchg \alpha\.\beta + \alpha\.\gamma$.
Staying with the notation introduced in Definition~\ref{umbraldot},
we would need to show that
\[ q_{\beta+\gamma,k}(\alpha) \eqv
     \sum_i {k\choose i} q_{\beta,i}(\alpha)q_{\gamma,k-i}(\alpha') \]
where $\alpha'\exchg \alpha$. This is not an equality. It only worked
when $\alpha\exchg a$ for $a\in\F$  because the substitution lemma
tells that 
$q_{\beta+\gamma,k}(\alpha) \eqv q_{\beta+\gamma,k}(a) $
and $\sum_i {k\choose i} q_{\beta,i}(\alpha)q_{\gamma,k-i}(\alpha') \eqv
     \sum_i {k\choose i} q_{\beta,i}(a)q_{\gamma,k-i}(a) $
which {\em is} equal to $q_{\beta+\gamma,k}(a)$.

The special case of $\alpha\.\gamma$ 
where $\alpha\exchg -n$ where $n$ is a positive integer is
Ray's definition in~\cite{Ra2} of ``negative umbral integers.''
Let $n$ be a positive integer.
Since $-n\.\gamma + n\.\gamma \exchg (-n+n)\.\gamma \exchg \varepsilon$,
the umbra $-n\.\gamma$ defined as above is precisely the same as the
umbral $-n\.\gamma$ defined in~\cite{RT2}.

The same techniques used in the preceding propositions proves the following.
\begin{proposition}
Let $a,c$ be in $\F$. If $p\in\F[\A]$ then
$a\.(cp) \exchg c(a\.p)$.
\end{proposition}
This points out that $-1\.\alpha$ is {\em not} in general exchangeable
with $-\alpha$. The latter is exchangeable with $-1(1\.\alpha)$.

Our definition of $p\.q$ does not allow for $p$ or $q$ to contain auxiliary
umbrae. Nevertheless, we would like to be able to manipulate expressions
that resemble $\alpha\.(\beta\.\gamma)$. Before we extend the notion of
an auxiliary umbra to handle this kind of construct, we prove the following
associativity result.
\begin{proposition}
Let $\alpha,\beta,\gamma$ be umbrae.
Define an umbra $\rho$ by $\rho\exchg\alpha\.\beta$ and
an umbra $\sigma$ by $\sigma\exchg\beta\.\gamma$. We have
$\rho\.\gamma \exchg \alpha\.\sigma$. 
\end{proposition}
Before presenting the proof, which is a quick calculation
with generating functions, we interpret the result probabilistically.
$A\.(B\.C)$ can be viewed as finding $A$, then running $A$ trials of 
$B\.C$, i.e. $A$ times we run a trial of $B$ and following each trial of 
$B$ we run that many trials of $C$. Then we add up all the trials of $C$.
In this interpretation $(A\.B)\.C$ differs only in that we run $A$ trials
of $B$ and then run all the trials of $C$ at once. 
We could extend this umbrae $\alpha,\beta,\gamma$ by viewing each side
as a polynomial in the  variables $\E[\alpha^i],\E[\beta^i],\E[\gamma^i]$.

Alternately, we argue as follows.
\begin{proof}
It suffices to check that
$e^{\rho\.\gamma\, z} \eqv e^{\alpha\.\sigma\, z}$.
If $a(z),b(z),c(z)\in\k[[z]]$ are given by 
$a(z)\eqv e^{\alpha\, z}$,
$b(z)\eqv e^{\beta\, z}$, and
$c(z)\eqv e^{\gamma\, z}$, then
this amounts to observing that each side is umbrally equivalent to
the composition $a(z)\comp\log(z)\comp b\comp \log(z)\comp c(z)$.
\end{proof}

With this lemma in hand, the following definition makes sense.
\begin{definition}
Given umbral polynomials $p_1,\ldots,p_n\in\k[\A_0]$, inductively
define the auxiliary umbra $p_1\.p_2\.\ldots\.p_n\in \A_1$ by
$p_1\.p_2\.\ldots\.p_n \exchg p_1\.\rho$ where
$\rho\exchg p_2\.\ldots\.p_n$.
\end{definition}

\section{Presentations for sequences of binomial type}

\subsection{Sequences of binomial type and sums of umbrae}

The notion of a sequence of binomial type is a direct generalization
of equation~\ref{bintake1exp}.
\begin{definition}
A sequence of polynomials $p_0(x),p_1(x),p_2(x),\ldots$ with $p_n(x)$ of degree $n$
is of {\em binomial type} when it satisfies
\begin{equation}\label{bineqn}
 p_k(x+y) = \sum_i {k\choose i} p_i(x)\, p_{k-i}(y) \, .
\end{equation}

Such a sequence is {\em normalized} when $p_1(x)$ is monic (equivalently $p_1(x)=x$).
\end{definition}

Equation~\ref{bintake1exp} arose directly as the umbral expansion 
of the identity (equation~\ref{bintake1}) 
that $(m+n)\.\alpha \exchg m\.\alpha + n\.\alpha'$. Recall that, for
any umbra $\gamma$ and any element $x\in\F$, $E[(x\.\gamma)^n]$ 
is a polynomial in $x$ and that 
$(x+y)\.\gamma \exchg x\.\gamma + y\.\gamma$
where $y$ is also in $\F$.
As in the proof of proposition~\ref{leftLin},
raising both sides of the preceding equality to the $n$th power and applying $\E$
implies the if direction of the following.
\begin{theorem}\label{probBin}
Let $p_0(x),p_1(x),p_2(x),\ldots$ be a sequence of polynomials
where $p_n(x)$ has degree $n$. This is as sequence of binomial
type iff it is umbrally represented by $x.\gamma$ for some umbra
$\gamma$.
\end{theorem}
\begin{proof}
By the remarks preceding the theorem, it suffices to show that any
sequence of binomial type can be so represented.
By standard results, which are briefly sketched below, it suffices to
show that choosing $\gamma$ appropriately allows us to choose the sequence
$D_x E[(x\.\gamma)^1]\evalAt_{x=0}\,, D_x E[(x\.\gamma)^2]\evalAt_{x=0}\,,
D_x E[(x\.\gamma)^3]\evalAt_{x=0}\,,\,\ldots$ arbitrarily.

It is enough to observe that 
equation~\ref{ProbCoeffs} tells us that the coefficient of $x$ in
$(x\.\gamma)^k$ is $\gamma^k + R$ where $R$ is depends only on $k$ and
$\E[\gamma],\ldots,\E[\gamma^{k-1}]$.
\end{proof}

For completeness, we sketch the fact that knowing $p'_1(0),p_2'(0),\ldots$
determines a sequence $p_0(x),p_1(x),p_2(x),\ldots$ of binomial type.
Replacing $y$ with $0$ in equation~\ref{bineqn} and recalling that,
by degree considerations the $p_i(x)$ are linearly independent, which tells us that
$p_0(0)=1$ and that $p_i(0)=0$ for $i>0$. Taking the
derivative of equation~\ref{bineqn} with respect to $y$ and and setting $y$ to $0$
gives $p_k'(x)=\sum_{i=0}^{k-1} {k\choose i} p_i(x) p_{k-i}'(0)$. Since 
$p_k(0)=\delta_{k,0}$ this determines $p_k(x)$.

Following~\cite{RKO}, the {\em umbral composition}
$a({\bf b}(x))$,
of two polynomial sequences $a_0(x),a_1(x),\ldots$ and $b_0(x),b_1(x),\ldots$
is the the sequence $T(a(0)),T(a(1)),\ldots$, where $T:\k[x]\rightarrow\k[x]$
is the linear operator defined by $T(x^i)=b_i(x)$ for all $i$.
An {\em umbral operator} is defined  to be a linear
operator $U:\k[x]\rightarrow\k[x]$ such that the sequence
$U(1),U(x),U(x^2),\ldots$ is of binomial type. The following
corollaries are immediate.
\begin{corollary}
A linear operator $U:\k[x]\rightarrow\k[x]$
is an umbral operator iff there exists an umbra $\gamma$ such
that $U\big(r(x)\big)=\E[\big(r(x)\big)\.\gamma]$
for all $r(x)\in\Q[x]$. \qed
\end{corollary}

\begin{corollary}
Let $p_0,p_1,\ldots$ and $q_0,q_1,\ldots$ be sequences of binomial type
represented by $x\.\alpha$ and $x\.\beta$ respectively. The umbral composition
of these sequences, $p_0({\bf q}), p_1({\bf q}),\ldots$ is represented by
$x\.\beta\.\alpha$. \qed
\end{corollary}
This makes obvious the fact from~\cite{RKO} that the umbral composition of
two sequences of binomial type is also of binomial type.

\subsection{Generalized Abel polynomials}\label{abelsect}

One of the best known sequences of binomial type has as its degree~$n$
polynomial the {\em Abel polynomial} $x(x+na)^{n-1}$ where
$a$ is a constant. Generalizing $a$ to be an arbitrary umbra $\alpha$ and replacing
$na$ with $n\.\alpha$ yields the following.

\begin{theorem}[\cite{RST}]\label{abel} 
Let $p_n(x)\in \k[x]$ be a sequence of polynomials with $p_1(x)=x$ and $p_n(x)$ of degree $n$.

The sequence $p_n(x)$ is of binomial type
iff there exists an umbra $\alpha$ such that
\[ p_n(x) \eqv  x(x+n\.\alpha)^{n-1} \,.\]
\end{theorem}
The proof in~\cite{RST} closely parallels the proof that the original Abel
polynomials are of binomial type. Here we provide a combinatorial proof.
\begin{proof}
To start with, assume that $\E[\alpha^i]$ is always an integer and that
$\alpha_1,\ldots,\alpha_n$ are distinct umbrae all exchangeable with $\alpha$.
We start by interpreting $(x+\alpha_1+\cdots+\alpha_n)^{n-1}$ as a generating
function for sequences of length $n-1$ on $n+1$ symbols. By the Pr\"ufer correspondence
(see for example~\cite{St})
this is a generating function for the number of labeled free trees on $n+1$ vertices where each
tree is counted with weight $x^{d_0}\prod_l \alpha_l^{d_l}$ where vertex $l$ has degree $d_l+1$.
So $x(x+\alpha_1+\cdots+\alpha_n)^{n-1}$ is the generating function for labeled rooted
trees on $n+1$ vertices where the same weight indicates that vertex $l$ has outdegree $d_l$.
This says that $\E[x(x+\alpha_1+\cdots+\alpha_n)^{n-1}]$ is the generating function where the
coefficient of $x^k$ counts the number of labeled trees on $n+1$ vertices where the root
has degree $k$ and each non-root vertex with outdegree~$i$ can be colored in any of $a_i$ ways.
Equivalently, $\E[x(x+\alpha_1+\cdots+\alpha_n)^{n-1}]$ counts the number of planted forests on 
$n$ vertices where each vertex with outdegree~$i$ can be colored in any of $a_i$ ways and where
each tree in the forest can itself be colored in any of $x$ ways. Let's call this structure a
$(x,\alpha)$ degree-colored forest on $n$ vertices.

So counting the number of ways to form a $(x+y,\alpha)$ degree-colored forest on $n$ vertex 
by the number of vertices, $i$, in the trees which were colored in one of the first $x$ ways  gives
\[ (x+y)(x+y+n\.\alpha)^{n-1} \eqv
\sum_i {n\choose i} (x)(x+i\.\alpha)^{i-1}\cdot (y)(y+{n-i}\.\alpha)^{n-i-1} \, .
\]
This fact for all positive integers $x,y,a_1,a_2,\ldots$ implies 
equation~\ref{bineqn} as a polynomial identity.

To see that indeed all normalized sequences of binomial type arise in 
this fashion, it suffices, by the remarks after Theorem~\ref{probBin}, to observe that the
sequence $p_2'(0),p_3'(0),\ldots$ can be chosen arbitrarily.
Indeed $p_n'(0)\eqv(n\.\alpha)^{n-1}\eqv n\alpha^{n-1}+R$
where $R$ is a linear combination of $\alpha_1,\ldots,\alpha_{n-1}$.
\end{proof}
The interpretation of $x(x+n\.\alpha)^{n-1}$ as a generating function for
colored forests was suggested to the author by Nigel Ray. It 
generalizes the notion of reluctant functions developed by
Mullin and Rota in~\cite{MR} and is closely related to the 
chromatic polynomials in~\cite{RaSW}.
 
The calculations used in~\cite{RST}  to prove Theorem~\ref{abel} show the following. 
\begin{proposition}
For any umbra $\alpha$ and any $n\ge1$ we find that  
\[
\displaylines{\hfill D_x e^{-1\.\alpha\,D_x} 
\left(x(x+n\.\alpha)^{n-1}\right) \eqv nx(x+(n-1)\.\alpha)^{n-2} \, .
\hfill\qedbox}\]
\end{proposition}

\begin{corollary}
A sequence $p_0(x),p_1(x),\ldots$ of polynomials, $p_n(x)$ of degree $n$ is a sequence of
binomial type iff there exists a formal power series $g(t)\in\k[[t]]$ with $g(0)=0$ and
$g'(0)\ne0$ such that
$g(D_x)\big(p_n(x) \big)=np_{n-1}(x)$ for all $n\ge1$.
\end{corollary}
\begin{proof}
If $p_1(x)$ is monic or $g'(0)=1$, the result follows immediately from the preceding proposition
and the remarks after Theorem~\ref{probBin}.

If the sequence is of binomial type and $p_1(x)=ax$, then so is the sequence whose
$n$th term is $p_n(x)/a^n$; if $g(t)$ is the series associated to this new sequence,
then $g(t)/a$ works for the original. The converse follows similarly.
\end{proof}
Formal power series of the sort described above are called {\em delta operators}
and the correspondence between them and their associated sequences of binomial type
was established in~\cite{RKO}.

 The transfer formula also arises as an immediate corollary. Since we now know that
the operator $Q\eqv De^{-1\.\alpha D_x}$ is associated to the sequence
presented as $p_n(x)\eqv x(x+\alpha_1+\cdots+\alpha_n)$, it follows that
$p_n(x)\eqv x e^{(\alpha_1+\cdots+\alpha_n)D_x} \left(x^{n-1}\right)=
x e^{\alpha_1 D_x}\cdots e^{\alpha_n D_x} \left(x^{n-1}\right) \eqv 
x\left({Q\over D}\right)^{-n}x^{n-1}.$

With the definition of a delta operator in hand, we recall the first expansion theorem
from~\cite{RKO}. If $Q$ is a delta operator associated to a sequence $p_0(x),p_1(x),\ldots$
of binomial type and if $T=f(D_x)$ for some formal power series
$f(t)\in\k[[t]]$, then
\begin{equation}
T=\sum_{k\ge0} Tp_k(0) {Q^k\over k!} \, .
\end{equation}
The usual proof employs the binomial expansion for $p_n(x+y)$ to verify the identity
$f(D_y)=\sum_{k\ge0} f(D_y)p_k(y) {Q^k\over k!}$
on $p_n(x+y)$ for all $n\ge0$. Setting $y=0$ gives the desired result.
The fact that $f(D_y)$ and the shift operator $e^{xD_y}$ commute is used freely.

We have recalled the expansion theorem in order to derive the identity
\begin{equation}\label{invseries}
D_x = \sum_{k\ge0} p_k'(0) {Q^k\over k!}
\end{equation}
where $p_0(x),p_1(x),\ldots$ is the sequence of binomial type associated to 
a delta operator, $Q$. Thus if $Q=f(D)$, we recover the fact that
$p_k'(0)$ is the coefficient of $t^k/k!$ in $f^{\langle -1\rangle}(t)$
the power series inverse, under composition, to $f(t)$. This relationship
together with the transfer formula is used to prove the Lagrange inversion
formula. See~\cite{MR,RST} for such derivations.

\subsection{Generalized rising factorials}\label{risfactsect}

Our next presentation generalizes the binomial type sequence of rising factorials,
$x(x+1)\ldots(x+n-1)$. More generally, it is well known that 
the sequence $p_n(x)=x(x+a)\ldots(x+(n-1)a)$ is of binomial type for all constants 
and that its associated delta operator for $a\ne0$ is $(I-e^{-aD})/a$. For the rising
factorials, this is the backwards difference operator $f(x)\mapsto f(x)-f(x-1)$.
Our proofs closely follow those for the usual rising factorials.

\begin{theorem}\label{ris.fact}
Let $p_n(x)\in \k[x]$ be a sequence of polynomials with $p_1(x)=x$ and $p_n(x)$ of degree $n$.

The sequence $p_n(x)$ is of binomial type iff
there exists an umbra $\mu$ such that
\[ p_n(x) \eqv x(x+\mu_1)(x+\mu_1+\mu_2)(x+\mu_1+\mu_2+\mu_3)\cdots(x+\mu_1+\ldots+\mu_{n-1}) \]
where $\mu_1,\ldots,\mu_{n-1}$ are distinct umbrae exchangeable with $\mu$.
\end{theorem}

\begin{proof}
We start by showing that all such presentations are of binomial type.

By induction we have that
\begin{eqnarray*}
\lefteqn{(x+y)(x+y+\mu_1)(x+y+\mu_1+\mu_2)\cdots(x+y+\mu_1+\ldots+\mu_{n-1}) \eqv }\\
& & y(y+x+\mu)(y+x+\mu+\mu_1)(y+x+\mu+\mu_1+\mu_2)\cdots
(y+x+\mu+\mu_1+\ldots+\mu_{n-2}) +\\
& & \quad +
x(x+y+\mu)(x+y+\mu+\mu_1)(x+y+\mu+\mu_1+\mu_2)\cdots
(x+y+\mu+\mu_1+\ldots+\mu_{n-2})
\\
&\eqv&
yp_{n-1}\Big(x+(y+\mu)\Big) + xp_{n-1}\Big((x+\mu)+y\Big)\\
&=& 
y\sum_i {n-1\choose i} p_{n-1-i}(x) p_i(y+\mu) \quad+\quad
  x\sum_i {n-1\choose i} p_{n-1-i}(x+\mu) p_i(y) \\
&\eqv&
\sum_i {n-1\choose i} p_{n-1-i}(x) 
      y(y+\mu)(y+\mu+\mu_1)\cdots(y+\mu+\mu_1+\ldots+\mu_{i-1}) 
\quad+\\ & &\quad+
  \sum_i {n-1\choose i} x(x+\mu)(x+\mu+\mu_1)\cdots(x+\mu+\mu_1+\ldots+\mu_{n-2-i}) 
      p_i(y) \\
&\eqv&
\sum_i {n-1\choose i} p_{n-1-i}(x) 
      p_{i+1}(y)
\quad+\quad
  \sum_i {n-1\choose i} p_{n-i}(x)
      p_i(y) \\
&\eqv&
\sum_{i=1}^n {n-1\choose i-1} p_{n-i}(x)  p_{i}(y)
\quad+\quad
  \sum_i {n-1\choose i} p_{n-i}(x)  p_i(y) \\
&\eqv&
  \sum_i {n\choose i} p_{n-i}(x)  p_i(y) \,. \\
\end{eqnarray*}

As before, all normalized sequences of binomial type arise in 
this fashion, since the
sequence $p_2'(0),p_3'(0),\ldots$ can be chosen arbitrarily.
Observe that $p_n'(0)\eqv\mu_1(\mu_1+\mu_2)\cdots(\mu_1+\cdots+\mu_{n-1})\eqv\mu_1^{n-1}+R$
where each term of $R$  has degree less than $n-1$.
\end{proof}

\begin{proposition}
Let $\mu$ be an umbra and let $\mu_1,\mu_2,\ldots$ be distinct umbrae exchangeable with $\mu$.
Consider the sequence of binomial type $p_n(x)\eqv x(x+\mu)\cdots(x+\mu_1+\cdots+\mu_{n-1})$
presented by the generalized rising factorials. The corresponding delta operator is
$q(D_x)$ where $q(t)\in\k[[t]]$ is given by $q(t)\eqv {e^{-1\.\mu\, D_x} -I\over {-1\.\mu}}$. 
\end{proposition}
\begin{proof}
Evaluating $q(D_x)$ on $p_n(x)$, we find
\begin{eqnarray*}
{e^{-1\.\mu\, D_x}\over -1\.\mu} p_n(x) &=& {p_n(x+-1\.\mu)-p_n(x)\over-1\.\mu}\\
&=& \sum_{i\ge1}{n\choose i} p_{n-i}(x) {p_i(-1\.\mu)\over -1\.\mu}\\
&\eqv& \sum_{i\ge1}{n\choose i} p_{n-i}(x)(-1\.\mu+\mu_1)(-1\.\mu+\mu_1+\mu_2)\cdots
                   (-1\.\mu+\mu_1+\cdots+\mu_{i-1})\\
&\eqv& {n\choose 1}p_{n-1}(x)\cdot1 \,.
\end{eqnarray*}
Here the last line follows from the substitution lemma by replacing $-1\.\mu+\mu_1$ with 
$\varepsilon\exchg0$.
\end{proof}

Just as the transfer formula arose from the generalized Abel presentation, the
Rodrigues formula arises from the presentation by generalized rising factorials.
Preserve the notation from the preceding proposition. We have
\begin{eqnarray*}
p_n(x)&\eqv& x(x+\mu_1)+\cdots(x+\mu_1+\cdots+\mu_{n-1})\\
&\eqv& x\cdot e^{\mu\, D_x} \big(x(x+\mu_2)\cdots(x+\mu_2+\cdots+\mu_{n-2})\big)\\
&\eqv& x\cdot e^{\mu\, D_x} p_{n-1}(x)\\
&\eqv& x\cdot \big(q'(D)\big)^{-1} p_{n-1}(x)\, ,\\
\end{eqnarray*}
where the last line follows since $q'(t)\eqv e^{-1\.\mu\, t}$.

We close our consideration of generalized rising factorials by observing that
the presentation $p_n(x)\eqv x(x+\mu_1)\cdots(x+\mu_1+\cdots+\mu_{n-1})$
immediately suggests a combinatorial interpretation along the lines of reluctant
functions. 
\begin{proposition}\label{colorFactorials}
Let $x$ be a nonnegative integer and let 
$\mu$ be an umbra such that $\mu^i\eqv m_i$ where each $m_i$ is a nonnegative integer.
The value of $p_n(x)\eqv x(x+\mu_1)\cdots(x+\mu_1+\cdots+\mu_{n-1})$ is
the number of labeled forests on $n$ vertices where each tree is assigned one of $x$ colors,
each vertex of outdegree~$j$ is assigned one of~$m_j$ colors, and where each parent vertex has
a smaller label than each of its children.
\end{proposition}
\begin{proof}
It suffices to observe that
we can construct such a tree by specifying a function mapping each vertex to its parent
and then choosing colors. Choosing $\mu_j$ from the $i$th multiplicand corresponds to
requiring the function to map vertex~$i$ to vertex~$j$. Choosing $x$
of course indicates that vertex $i$ is a root.
\end{proof}

\section{Presentations for Sheffer sequences.}

Recall that a sequence of polynomials $s_0(x),s_1(x),s_2(x),\ldots$ with
$s_n(x)$ of degree $n$ is said to be a {\em Sheffer sequence} with respect
to a delta operator $Q$ (or with respect to the associated sequence of binomial type)
when $Q s_n(x)=n s_{n-1}(x)$ for all $n\ge0$. We will
call a Sheffer sequence {\em normalized} when $s_1(x)$ is monic. 

Let $Q=f(D_x)$ be a delta operator. Let $p_0(x),p_1(x),\ldots$ be the associated sequence
of binomial type. Since $Q$ and $e^{\beta D_x}$ commute, we have
$Q p_n(x+\beta)\eqv n p_{n-1}(x+\beta)$ thus $p_n(x+\beta)$ is a Sheffer sequence for $Q$.
Since $Q/D$ is invertible, $s_n(0)$ and $Qs_{n}(x)$ determine $s_n(x)$. Hence  
$s_1(0),s_2(0),\ldots$ determines any Sheffer sequence with respect to $Q$.
But because the the $p_i(x)$ have different degrees, any such sequence arises from 
the umbral presentation $s_n(x)\eqv p_n(x+\beta)$ for suitable choice of umbra $\beta$.
We have proved the following.
\begin{proposition}
If $p_0(x),p_1(x),\ldots$ is a sequence of binomial type, then all associated
Sheffer sequences are umbrally presented by the sequence
 $p_0(x+\beta),\, p_1(x+\beta),\, p_2(x+\beta),\,\ldots$
for some $a\in\k$ and some umbra $\beta$.\qed
\end{proposition}

The following is immediate from the preceding and the presentation results 
in the preceding section.
\begin{corollary}
Let $s_0(x),s_1(x),\ldots$ be a sequence of polynomials in $\k[x]$
and let $p_0(x),p_1(x),\ldots$ in $\k[x]$ be a sequence of binomial type.
The following are equivalent:
\begin{enumerate}
\item The sequence $s_0(x),s_1(x),\ldots$ is Sheffer 
with respect to $p_0(x),p_1(x),\ldots$.

\item \label{SheffRV}
There exist umbrae $\beta,\gamma$ such that 
$p_n(x)\eqv (x\.\gamma)^n$
and 
$s_n(x)\eqv \big((x+\beta)\.\gamma\big)^n$.

\item There exist umbrae $\beta,\alpha$ such that 
$p_n(x)\eqv x(x+n\.\alpha)^{n-1}$
and 
$s_n(x)\eqv (x+\beta)(x+\beta+n\.\alpha)^{n-1}$.

\item There exist umbrae $\beta,\mu$ such that 
$p_n(x)\eqv x(x+\mu_1)\cdots(x+\mu_1+\cdots+\mu_{n-1})$
and 
$s_n(x)\eqv (x+\beta)(x+\beta+\mu_1)\cdots(x+\beta+\mu_1+\cdots+\mu_{n-1})$.

\end{enumerate}
\qed
\end{corollary}

Since $(x+y+\beta)\.\gamma \exchg x\.\gamma + (y+\beta\.\gamma)$,
the standard expansion result for Sheffer sequences,
\[ s_n(x+y)=\sum_{i=1}^n {n\choose i} p_i(x)s_{n-i}(y) \,,\]
follows immediately from part~\ref{SheffRV} of this corollary.

\section{Multiplicative sequences} 

We present our final results as easy applications of the preceding constructions.
We have relied on the umbral relation $(x+y)\.\gamma \exchg x\.\gamma + y\.\gamma$
to find a presentation for sequences of binomial type. It is natural to ask what
happens if we replace $x$ and $y$ themselves by umbrae. Fix umbrae $\alpha,\beta$ such
that $\alpha^n\eqv a_n$ and $\beta^n\eqv b_n$. We will consider sequences of polynomials
in multiple variables $a_1,a_2,\ldots$ and $b_1,b_2,\ldots$.

Fix an umbra $\gamma$ and
define $K_m(a_1,a_2,\ldots)\eqv (\alpha\.\gamma)^m$.
By equation~\ref{ProbCoeffs} $K_m(a_1,a_2,\ldots)$ has degree~$m$, is linear in the variables $a_i$
and only depends on the variables $a_1,\ldots,a_m$.

>From the point of view of generating functions, the fairly simple umbral relation
$(\alpha+\beta)\.gamma\exchg \alpha\.\gamma+\beta\.\gamma$ now becomes:
\begin{proposition}
Define  polynomials $K_m\in\k[a_1,\ldots,a_m]$ by $K_m(a_1,a_2,\ldots,a_m)\eqv (\alpha\.\gamma)^m$.
If \[ \sum_{k\ge0} c_k {z^k\over k!} = \left( \sum_{i\ge0} a_i {z^i\over i!} \right)
\left( \sum_{j\ge0} b_j {z^j\over j!} \right) \, ,\]
then 
\[ \displaylines{ \hfill \sum_{k\ge0} K_k(c_1,\ldots,c_k) {z^k\over k!} = 
\left( \sum_{i\ge0} K_i(a_1,\ldots,a_i) {z^i\over i!} \right)
\left( \sum_{j\ge0} K_j(b_1,\ldots,b_j) {z^j\over j!} \right) \, .
\hfill\qed}\]
\end{proposition}

We generalize this ``nice'' behavior of generating functions under multiplication
as follows. Suppose that $K_0,K_1,\ldots$ is a sequence of polynomials in
the variables $t_1,t_2,\ldots$. For $r\in\k[\A]$, denote by $K_m(r)$ the
evaluation, $K_m(\E[r],\E[r^2],\E[r^3],\ldots)$, of the polynomial $K_m$.
Define the sequence $K_0,K_1,\ldots$
to be {\em \multiplicative} if whenever $\alpha,\beta,\chi$ are umbrae such that
$\alpha + \beta \exchg \chi$, we have $\rho+\sigma\exchg \tau$ where
$\rho^m\eqv K_m(\alpha)$, $\sigma^m\eqv K_m(\beta)$, and $\tau^m\eqv K_m(\gamma)$. 

Since $n\.(\alpha+\beta)\exchg n\.\alpha + n\.\beta$,
the sequence of polynomials $K_m$ in $a_1,a_2,\ldots$ defined by 
$K_m(a_1,\ldots,a_m)\eqv (n\.\alpha)^m$ is \multiplicative.

These constructions can be generalized as follows.
\begin{proposition}
Let $\alpha$ be an umbra with $\alpha^i\eqv a_i$.
Let $l$ be a positive integer and $\gamma$ and umbra. For $c_1,\ldots,c_l\in\k$,
define the polynomial $K_m(a_1,\ldots,a_m)$ in $\k[a_1,\ldots,a_m]$ by
$K_m(a_1,\ldots,a_m)\eqv (\sum_{i=1}^l i\.(c_i\alpha)\.\gamma)^m$.
The sequence $1,K_1,K_2,\ldots$ is \multiplicative and $K_m$ has total
degree $m$ when $a_i$ is given degree $i$.

If $\gamma\exchg 1$, then $K_m$ is homogeneous in the above grading.\qed
\end{proposition}

If $K_0,K_1,\ldots$ is \multiplicative and homogeneous, then
$L_0,\, L_1,\, \ldots$ where $L_i=K_i/i!$ is a 
$m$-sequence or multiplicative sequences in the sense of Hirzebruch~\cite{Hirz},
namely if 
 \[ \sum_{k\ge0} c_k {z^k} = \left( \sum_{i\ge0} a_i {z^i} \right)
\left( \sum_{j\ge0} b_j {z^j} \right) \, ,\]
then 
\[  \hfill \sum_{k\ge0} L_k(c_1,c_2,\ldots,c_k) {z^k} = 
\left( \sum_{i\ge0} L_i(a_1,a_2,\ldots,a_i) {z^i} \right)
\left( \sum_{j\ge0} L_j(b_1,b_2,\ldots,b_j) {z^j} \right) \, .
\]

\section{Open problems}

We close with a brief list of open problems and areas for future work.
\begin{enumerate}
\item Generalize Theorems~\ref{abel} and~\ref{ris.fact} and their corollaries
by finding umbral presentations   
corresponding to other well-known sequences of binomial type. 

\item Determine which sequences of binomial type over the integers
and which sequences of integral type (see~\cite{BBN}) may be presented
in the above fashions. Find general presentation formulae for these
situations.

\item Find umbral presentation theorems for the various generalizations of the umbral
calculus (see for instance~\cite{LR,Loeb}).

\item Give conditions for a sequence of binomial type to be presentable
as $x\.G$ where $G$ is a random variable rather than an arbitrary umbra.

\end{enumerate}

\section{Acknowledgments} 
I am greatly indebted to Gian-Carlo Rota who first suggested that
the formalization of umbrae developed in~\cite{RT1,RT2} was amenable 
to interpretation as a generalization of random variables. I am likewise
indebted to Nigel Ray; the influence of my conversations with him
can be traced throughout subsections~\ref{abelsect} and~\ref{risfactsect}.

\end{document}